\documentclass[aop,MSNbibl,citesort,dvips]{arximspdf}

\doi{10.1214/10-AOP593}
\volume{39}
\issue{2}
\pubyear{2011}
\firstpage{407}
\lastpage{416}

\begin{document}
\begin{frontmatter}

\title{T. E. Harris' contributions to interacting particle systems and
percolation}
\runtitle{Harris IPS and percolation}

\begin{aug}
\author[A]{\fnms{Thomas M.} \snm{Liggett}\corref{}\thanksref{t1}\ead[label=e1]{tml@math.ucla.edu}}
\runauthor{T. M. Liggett}
\affiliation{University of California at Los Angeles}
\address[A]{Department of Mathematics\\
University of California\\
520 Portola Plaza\\
Los Angeles, California 90095\\
USA\\
\printead{e1}} 
\end{aug}

\thankstext{t1}{Supported in part by NSF Grant DMS-03-01795.}

\received{\smonth{8} \syear{2010}}

%
\begin{abstract}
Interacting particle systems and percolation have been among the most
active areas of probability theory over the past half century. Ted
Harris played an important role in the early development of both
fields. This paper is a bird's eye view of his work in these fields,
and of its impact on later research in probability theory and
mathematical physics.
\end{abstract}

\setattribute{keyword}{AMS}{AMS 2000 subject classification.}

\end{frontmatter}

\section{Introduction}
Ted's passing was a great loss to me personally and professionally, as
well as to probability theory in Southern California and beyond. For
three decades, he and I were the primary probabilists at USC and UCLA,
respectively. Our approach to mathematics and our mathematical
interests were very similar. My wife, Chris, and I enjoyed many
wonderful social occasions at the Harris home in Beverly Hills with Ted
and Connie. They were great hosts. The conversation was always
stimulating. They were interested in so many things!

Ted initiated the Southern California Probability Symposium in about
1970. It is probably the oldest meeting of its type in the US, and
continues to provide an exciting place for interactions among Southern
California probabilists to this day. In two of the years, 1989 and
2007, the meeting was dedicated to Ted. The first was on the occasion
of his 70th birthday, and the second was in his memory. Even after his
retirement, Ted attended the USC probability seminar
regularly---always sitting in the front row, and asking perceptive questions.

Few mathematicians have had a greater ratio of number of ideas to
number of papers. Ted wrote fewer papers (about 30) than many prominent
mathematicians, but each is a jewel. In each of the several areas in
which he worked, he was among the first in the field. He had an uncanny
ability to sense which problems would lead to major developments. His
taste was impeccable. After having a significant impact on one area, he
would go on to another topic, leaving it to others to flesh out the
subject. His work went a long way toward shaping the growth of
probability theory in the second half of the twentieth century. I hope
to do justice to his many ideas on interacting particle systems and
percolation in this brief article. I will mention only a few later
papers in which his influence can be seen---there are many others.

\section{Percolation} In the standard percolation model, bonds in
$Z^d$ are independently labeled open with probability $p$ and closed
with probability $1-p$. One then asks whether the subgraph of $Z^d$
obtained by retaining only the open bonds contains an infinite
connected component. The answer of course depends on the value of
$p$---it is yes if $p$ is sufficiently large, and no if $p$ is small. The
value at which the answer changes is known as the critical value, $p_c$.

The mathematical theory of percolation is generally viewed as beginning
in 1957 with paper \cite{BH} by Broadbent and Hammersley. Only three
years later, Ted proved in \cite{H60} that in two dimensions, $p_c\geq
\frac12$. In fact, he proved that there is no infinite cluster when
$p=\frac12$. Prior to that, the best result on the critical value was
$0.35\leq p_c\leq0.65$. The primary tools he used were a correlation
inequality (more on this in Section \ref{corr} below), the
self-duality of the two-dimensional lattice, and path intersection
arguments. It took another 20 years for Kesten to prove in \cite
{Ke} that $p_c=\frac12$, using among other techniques, more refined
path intersection arguments. It takes only a glance at Grimmett's book
\cite{Grim} or the more recent \cite{BR} by Bollob\'as and Riorden to
get a sense of how big and important percolation has become since Ted's
pioneering work. Recent work on SLE scaling limits for percolation and
other models is one of the most exciting developments in modern
probability theory---see \cite{BR,CN,La} and \cite
{Sm}, for example. In fact, one of the 2010 Fields Medals was awarded
to Smirnov for his work in the field. Again, path intersection
arguments play an important role. Ted did not return explicitly to
percolation after \cite{H60}, but percolation ideas were to play a
major role in his later work in interacting particle systems.

\section{The contact process} Contact processes constitute one of the
two or three major classes of interacting particle systems. They play
somewhat the same role in this area that Brownian motion plays in the
theory of stochastic processes in Euclidean space: they are simple to
describe, they have many of the useful properties that other systems in
the field may or may not have---in this case, self-duality,
attractiveness and additivity---and they lead to challenging
mathematical problems. Contact processes were first introduced and
studied by Ted Harris in \cite{H74}. Literally hundreds of papers have
been written about them in the past 35 years. Few mathematicians have
been credited with starting a field that would become as important as
this one.

The basic contact process on $S=Z^d$ is a Markov process $\eta_t$ on
$\{0,1\}^S$, which can be thought of as a model for the spread of
infection. A configuration $\eta\in\{0,1\}^S$ represents the state in
which certain sites are infected [those for which $\eta(x)=1$]; the
others are healthy. The value at a site $x\in S$ changes from 1 to 0 at
rate 1 (i.e., infected sites recover after a unit exponential time),
and from 0 to 1 at a rate proportional to the number of infected
neighbors. The constant of proportionality is $\lambda$. Of course,
the configuration $\eta\equiv0$ is a trap for the
process---infections cannot appear spontaneously. When there are only finitely
many infected sites, the state of the system is usually denoted by
$A_t=\{x\dvtx\eta_t(x)=1\}$.

A principal reason for interest in the contact process is that it, like
the percolation model, can have two different types of behavior,
depending on the value of a parameter---$\lambda$ in this case. It
can survive, in the sense that the survival probability is positive,
%
%
\begin{equation}\label{surv}\pi(A)\equiv P^A(A_t\neq\varnothing
\mbox{ for all }t)>0,\qquad A\neq\varnothing,
\end{equation}
even for initial configurations with finitely many infected sites, or
it can die out. In \cite{H74}, Ted did not talk in terms of the
critical value $\lambda_c$ that separates these two regimes, or even
note that this value was well defined. However, he did prove that the
process survives if $\lambda$ is large enough, so that $\lambda
_c<\infty$, and that $\lambda_c\geq\frac1{2d-1}$, with an
improvement to $\lambda_c\geq1.18$ in one dimension. By now, much
more is known: $1.53\leq\lambda_c<2$ in one dimension, and $\lambda
_c\sim1/2d$ in high dimensions. The actual value of $\lambda_c$ is
thought to be about 1.65 in one dimension, but nothing close to this is
rigorously known. (Most later results that are mentioned in this
article can be found in \cite{Li1} or \cite{Li2}.)

Ted was a fan of inequalities, as we will see in Section \ref{corr}.
In \cite{H74}, he proved some inequalities that are not of the
correlation type discussed there. It is fairly clear that the survival
probability $\pi(A)$ is increasing in $A$. It is less obvious that it
is submodular in the sense that
%
%
\begin{equation}\label{subm}\pi(A\cup B)+\pi(A\cap B)\leq\pi
(A)+\pi(B).
\end{equation}
This inequality played a role in his derivation of lower bounds for the
critical value.
Only much later was this result used in an essential way in \cite{LSS}
to compare a contact process with mutations to the basic contact
process, showing that the former process dies out whenever the latter
does. This comparison apparently cannot be carried out via more common
and intuitive coupling arguments. As far as I know, (\ref{subm}) was
not used in the intervening 34 years, even though it was generalized to
some extent, and has been used in some related contexts---see
\cite{MR}. Ted was again ahead of his time.

One of the most useful techniques in interacting particle systems is
duality, which expresses probabilities related to one process in terms
of probabilities related to another (dual) process. Forms of duality
had been used quite early in the study of Brownian motion and birth and
death chains. In the context of symmetric exclusion processes, duality
was discovered by Spitzer in \cite{Sp}, and has played an essential
role in that theory. In particular, it made possible a complete
description of the stationary distributions of the system. Such a
classification in the asymmetric case remains elusive.

In \cite{H76}, Ted looked at duality more generally, and discovered in
particular that the contact process is self-dual. (He used the word
``associate'' rather than ``dual.'') Self-duality for the contact process
is the identity
\[
P^\eta(\eta_t\equiv0\mbox{ on }A)=P^A(\eta\equiv0\mbox{ on
}A_t),\qquad |A|<\infty,
\]
which relates the contact process with infinitely many infections to
the process with finitely many infections. Letting $1$ denote the
configuration with all sites infected, this says in particular that
\[
P^1\bigl(\eta_t(x)=1\bigr)=P^{\{x\}}(A_t\neq\varnothing),
\]
so that survival in the sense of (\ref{surv}) is equivalent to
survival in the sense that
\[
\lim_{t\rightarrow\infty}P^1\bigl(\eta_t(x)=1\bigr)>0.
\]
When the process survives, the ``upper invariant measure'' $\nu$ is
defined as the limiting distribution as $t\rightarrow\infty$ of the
distribution at time $t$ of the system starting in configuration 1.
Thus, duality gives
\[
\nu\{\eta\dvtx\eta\not\equiv0\mbox{ on }A\}=\pi(A).
\]

Duality can be used to give a simple proof of the submodularity
property (\ref{subm}) that Ted had discovered earlier and proved by
coupling. It is obtained by integrating the elementary inequality
\[
1_{\{\eta\not\equiv0\ \mathrm{on}\ A\cup B\}}+1_{\{\eta\not\equiv
0\ \mathrm{on}\ A\cap B\}}\leq1_{\{\eta\not\equiv0\ \mathrm{on}\ A\}
}+1_{\{\eta\not\equiv0\ \mathrm{on}\ B\}}
\]
with respect to $\nu$. Of course, Ted's duality was not available to
him when he proved~(\ref{subm}).

There are a number of applications of duality in \cite{H76}, including
a proof of the fact that every translation invariant stationary
distribution for the contact process is a mixture of $\nu$ and the
point mass on $\eta\equiv0$. Now, we know that the translation
invariance assumption is not needed in this statement.

Another important technique in the field is known as the graphical
representation, or percolation substructure. The basic idea developed
in \cite{H78} is that it is very natural and useful to construct
processes like the contact process explicitly in terms of collections
of independent Poisson processes. There are many advantages to this
approach, including the possibility of constructing the process
starting from all potential initial configurations on the same
probability space. It also gives duality in an explicit way. In the
space--time graphical picture, the evolution of the dual process is seen
by reversing the time direction.

The graphical representation has played a crucial role in many proofs,
including the 1990 proof by Bezuidenhout and Grimmett \cite{BG} that
the critical contact process dies out. It is the underlying theme of
Griffeath's monograph \cite{Grif}, and is the basis of a lot of work
of Durrett and his collaborators on systems related to the contact
process---see his paper \cite{Du} in the volume dedicated to Ted's
70th birthday, for example.

In his paper, Ted proved a number of results using the graphical
representation. Here are two:

(a) Linear growth: for sufficiently large $\lambda$,
\[
P^A\biggl(\inf_{t>0}\frac{|A_t|}t>0 \Big| A_t\neq\varnothing
\mbox{ for all }t\biggr)=1.
\]
He points out that the $t$ in the denominator can probably be replaced
by the more plausible $t^d$.

(b) The individual ergodic theorem: for a large class of initial $\eta
$ and all continuous functions $f$,
\[
\lim_{T\rightarrow\infty}\frac1T\int_0^Tf(\eta_t)\,dt=\int
f\,d\nu\qquad\mbox{a.s.}
\]
Now we know that (a) holds (with $t^d$ in the denominator) for any
$\lambda>\lambda_c$. Statement~(b) has also been improved.

\section{Exclusion processes}

This represents another large part of the field of interacting particle
systems. An exclusion process is described via the transition
probabilities $p(x,y)$ for a discrete time Markov chain on a countable
set $S$. The process is again a continuous time Markov process on
$\{0,1\}^S$. This time, however, $\eta(x)=1$ means that site $x$ is
occupied; $\eta(x)=0$ means that it is vacant. A particle at $x$ waits
a unit exponential time, and then chooses a $y\in S$ with probabilities
$p(x,y)$. If $y$ is vacant, the particle moves there, while if $y$ is
occupied, it stays at~$x$. The process was introduced by Spitzer in
\cite{Sp}, and has been the subject of a very large number of papers in
both mathematics and physics---both rigorous and nonrigorous---since
then.

Again, Ted was in the game from the beginning. Prior to my general
existence theorem in \cite{Li0}, three mathematicians constructed
particle systems under various assumptions---Dobrushin \cite{Do},
Holley \cite{Ho} and Harris \cite{H72}. In his paper, Ted was
concerned with nearest neighbor exclusion processes on $Z^d$. He used
percolation ideas in his construction. Noting that it suffices to
construct a Markov process for an arbitrarily short interval of
time---the Markov property allows for an extension to all time---he showed
that for such short time periods, $Z^d$ breaks up into random finite
subsets that do not interact with one another during that time period.
On each of these finite subsets, the process is of course well defined.

A useful point of view in the study of symmetric [i.e., those
satisfying $p(x,y)=p(y,x)$] exclusion processes was pioneered by Ted,
and is known as ``stirring.'' This is closely related to the graphical
representation he introduced for contact-like processes in \cite{H78}.
The idea is that a Poisson process of rate $p(x,y)$ is associated with
each pair of sites $x,y$. At the event times of the Poisson process,
the ``contents'' of the two sites are exchanged. If both were empty or
both were occupied, nothing happens, since the particles are
indistinguishable. If exactly one site is occupied, the result is that
the particle at the occupied site moves to the other site. This again
constructs all exclusion process with arbitrary initial configurations
on the probability space of the Poisson processes. With this
construction, the system is realized as a collection of interacting
copies of the original Markov chain. Ted wrote about stirring in
\cite{H91}, and one of his students used it in \cite{Le}.

An important application of stirring occurred in a paper with another
connection to Ted's work. In \cite{H65}, he considered a system of
reflecting Brownian motions, one starting at each point of a unit
Poisson process on the line, and an extra one starting at the origin.
He defined reflection by saying that when two Brownian paths meet, they
interchange their paths, so that the particles maintain their original
ordering. He then proved that the position of the particle originally
at the origin satisfies a central limit theorem, but with scaling
$t^{1/4}$, rather than $t^{1/2}$. This led Spitzer in \cite
{Sp} to make a conjecture on the behavior of a ``tagged'' particle in a
symmetric exclusion process. Ted's USC colleague R. Arratia then proved
the conjecture in~\cite{Ar}, using stirring in an essential way. Here
is the result. Consider the exclusion process with $S=Z^1$ and
\[
p(x,x+1)=p(x,x-1)=\tfrac12.
\]
Initially there is a particle at the origin, and particles are placed
at other sites independently with probability $\rho\in(0,1)$. The
particle that started at the origin is the tagged particle. Its
position satisfies a central limit theorem, but again with the
nonstandard normalization. It turns out that this is (apparently) the
only case in which an unusual scaling occurs. Central limit theorems
for tagged particles in exclusion systems have been proved by Varadhan
and others with normalization $t^{1/2}$ in many cases in which
$p(x,y)$ is translation invariant on $Z^d$, including systems with mean
zero in \cite{KV} and \cite{Va} in any dimension (excluding Arratia's
case), and systems with nonzero mean in dimensions $d\geq3$ in \cite{SVY}.

\section{Correlation inequalities}\label{corr} Perhaps Ted's best
known and most influential result is the correlation inequality in
\cite{H60}. Amazingly, it is not even mentioned in the MathSciNet
review of that paper. Perhaps that is not so amazing after all. Who
would have known 50 years ago what an effect it would have?

To state it, let $S$ be a finite set, and consider the Bernoulli
measure $\nu_{\rho}$ on $\{0,1\}^S$ defined by
\[
\nu_\rho\{\eta\dvtx\eta(x)=1 \mbox{ for all }x\in T\}=\rho
^{|T|},\qquad T\subset S.
\]
A set $A\subset\{0,1\}^S$ is said to be increasing if $\eta\in A$ and
$\eta\leq\zeta$ imply that $\zeta\in A$. [$\eta\leq\zeta$ means
that $\eta(x)\leq\zeta(x)$ for all $x\in S$.] Ted's result is
%
%
\begin{equation}\label{corr1}A,B\mbox{ increasing implies }\nu_\rho
(A\cap B)\geq\nu_{\rho}(A)\nu_{\rho}(B).
\end{equation}
Actually, he only proved this for the increasing sets that arose in his
percolation problem, but that is a minor point.

Property (\ref{corr1}) is now usually stated in terms of increasing
functions rather than sets, and when applied to a general probability
measure, is called ``association.'' Thus, a probability measure $\mu$
on $\{0,1\}^S$ is said to be associated if
%
%
\begin{equation}\label{corr2}f, g\mbox{ increasing implies }\int
fg\,d\mu\geq\int f\,d\mu\int g\,d\mu.
\end{equation}
Ted's theorem then states that homogeneous product measures are
associated. Motivated by this result, as well as by a 1967 result of
Griffiths \cite{Griff}, Fortuin, Kasteleyn and Ginibre \cite{FKG}
proved a far reaching generalization that is known as the FKG theorem:
if the probability measure $\mu$ is strictly positive and satisfies
%
%
\begin{equation}\label{lattice}\mu(\eta\wedge\zeta)\mu(\eta\vee
\zeta)\geq\mu(\eta)\mu(\zeta),\qquad \eta,\zeta\in\{0,1\}^S,
\end{equation}
where $\eta\wedge\zeta$ and $\eta\vee\zeta$ denote the
coordinate-wise minimum and maximum of $\eta$ and~$\zeta$, respectively,
then $\mu$ is associated. In their paper, they mentioned the
understated nature of Ted's result, while recognizing its importance:
``While Harris' inequality seems to have drawn less attention than it
deserves$,\ldots.$'' Ted himself, with his usual modesty, said
``perhaps the
methods are also of some interest.'' Note that while (\ref{corr1}) is
far from obvious, the lattice condition (\ref{lattice}) is easy to
check for $\nu_\rho$. The FKG theorem has played a pivotal role in
the study of phase transitions in statistical physics over the past
four decades.

Ted's other paper on correlation inequalities \cite{H77} is short,
elegant, and has many consequences. It deals with implications among
the following three properties for a continuous time Markov process
$\eta_t$ on $\{0,1\}^S$:

(a) Preservation of association: if $\mu$ is associated, then so is
$\mu_t$, the distribution of $\eta_t$ with initial distribution $\mu$.

(b) All transitions are between comparable configurations.

(c) Attractiveness: if $f$ is increasing on $\{0,1\}^S$, then so is
$E^\eta f(\eta_t)$.

His theorem is that in the presence of (c), (a) and (b) are
equivalent. Much later, I proved in \cite{Li3} that if all transitions
are between configurations that differ at only one site, then (a)
implies (c).

An easy consequence of Ted's result is that the upper invariant measure
for the contact process $\nu$ is associated. To see this, note that
the point mass at $\eta\equiv1$ is associated, and the contact
process satisfies (b) and (c). Therefore, the distribution at time $t$
is associated. Now, let $t\rightarrow\infty$. The fact that $\nu$ is
associated is not a consequence of the FKG theorem, since it is known
that $\nu$ does not satisfy (\ref{lattice}).

Ted's theorem has been extended in the case of the contact process to
show that certain properties that lie between (\ref{lattice}) and
association are also preserved by the evolution---see Theorem 3.5 of
\cite{BHK} and Theorems 1.5 and 1.7 of \cite{Li3}. One consequence of
this is that the upper invariant measure $\nu$ percolates if $d\geq2$
and $\lambda$ is sufficiently large---see \cite{LS}.

It is interesting to note that (\ref{corr1}) follows from Ted's later
theorem: consider the spin system in which the coordinates $\eta_t(x)$
flip independently from 1 to 0 at rate $1-\rho$ and from 0 to 1 at
rate $\rho$. This process satisfies (b) and (c), and has limiting
distribution $\nu_\rho$ for any initial state. Therefore, $\nu_\rho
$ is associated.

Paper \cite{H77} has stimulated recent work on negative correlations
as well. The theory of negative correlations is more subtle than that
of positive correlations. One way to see this is that in the definition
of negative association, one cannot simply reverse the inequality in
(\ref{corr2}), as can be seen by taking $f=g$ there. One must add the
constraint that $f$ and $g$ depend on disjoint sets of coordinates.
With this definition, the negative version of the FKG theorem is false:
(\ref{lattice}) with the opposite inequality does not imply negative
association.

A possible version of Ted's 1977 theorem for negative association might
be that systems like the symmetric exclusion process preserve the
property of negative association. After all, if many particles are
known to be in one part of $S$, then fewer can be in other parts of
$S$. In hindsight, it is not too surprising that this is false. The
restriction to functions that depend on disjoint sets of coordinates in
the definition causes problems, since even if $f$ and $g$ satisfy this
constraint, $E^\eta f(\eta_t)$ and $E^\eta g(\eta_t)$ generally will
not. All is not lost, however.
In \cite{BBL}, a property that is stronger than negative association,
but still satisfied by product measures, is shown to be preserved by
symmetric exclusion processes.

\section*{Acknowledgments}
I appreciate comments from Rick Durrett and Geoffrey Grimmett on this paper.

%

%
\printaddresses

\end{document}